\documentclass[12pt]{amsart} 
\usepackage{amssymb}
\usepackage{fullpage}
\usepackage[colorlinks=true]{hyperref}

\newcommand{\bbQ}{\mathbb{Q}}

\newcommand{\bbZ}{\mathbb{Z}}

\newcommand{\cN}{\mathcal{N}}
\newcommand{\cO}{\mathcal{O}}

\DeclareMathOperator{\core}{sf}

\numberwithin{equation}{section}

\newtheorem{theorem}{Theorem}[section]
\newtheorem{conjecture}[theorem]{Conjecture}
\newtheorem{corollary}[theorem]{Corollary}
\newtheorem{lemma}[theorem]{Lemma}
\newtheorem{proposition}[theorem]{Proposition}

\theoremstyle{definition}

\newtheorem*{remark-nonum}{Remark}

\newcommand{\lemDio}{2.1}
\newcommand{\lemHypg}{2.2}

\newcommand{\eqDPrimeDef}{2.6} % definition of d'

\newcommand{\lemLB}{3.5}

\newcommand{\lemGap}{3.8}
\newcommand{\lemTheta}{3.9}

\newcommand{\lemGNd}{3.12}

\newcommand{\propRep}{3.1}

\newcommand{\subSectConstr}{2.1} % Construction of Approximations subsection
\newcommand{\subSectStepiii}{4.4} % r_0>1 and p/q \neq p_r0/q_r0

\newcommand{\eqQAbs}{3.25} % from proof of gap principle, used to get expression of |q|
\newcommand{\eqGPb}{3.28}   % gap principle part (b) proof lower bound for (1/y_i+1/y_j)
\newcommand{\eqOmegaApproxUB}{4.2} % upper bound for |w^{1/4}-approx|
\newcommand{\eqELB}{4.5} % lower bound for E
\newcommand{\eqQUB}{4.7} % upper bound for Q
\newcommand{\eqEllUB}{4.8}    % ell_0 upper bound in the prereqs
\newcommand{\eqYLUBUsingYK}{4.15} % upper bound for y_\ell for any r_0 and y_k

\newcommand{\corMainSqr}{1.4}

\newcommand{\subsectPrereqs}{4.1}
\newcommand{\eqYUBStepOne}{4.14}
\newcommand{\eqYUBStepThree}{4.22}
\newcommand{\eqUBStepFour}{4.28} % the upper bound for (stuff)<1

\newcommand{\corMainPP}{1.7}

\begin{document}

\title[Bounds on the number of squares \ldots]{Bounds on the number of squares in recurrence sequences: arbitrary $b$, III}

\author{Paul M Voutier}
\address{London, UK}
\email{Paul.Voutier@gmail.com}

\date{}

\begin{abstract}
We generalise our earlier work on the number of squares in binary recurrence
sequences, $\left\{ y_{k} \right\}_{k \geq -\infty}$. In the notation of our
previous papers, here we consider the case when $N_{\alpha}$ is any negative
integer and $y_{0}=b^{2}$ for any positive integer, $b$. We show that there
are at most $4$ distinct squares with $y_{k}$ sufficiently large. This allows
us to also show that there are at most $9$ distinct squares in such sequences
when $b=1,2$ or $3$, or once $d$ is sufficiently large.
\end{abstract}


\maketitle

\section{Introduction}
\label{sect:intro}

\subsection{Background}

In recent work \cite{V5,V6}, we developed a technique for bounding the number
of distinct squares in binary recurrence sequences. In these papers, we applied
our technique to sequences $\left( y_{k} \right)_{k \geq -\infty}$ with $y_{0}=1$
that arise from the solutions of generalised Pell equations, $X^{2}-dY^{2}=c$
with $c<0$. We were able to obtain best possible results for most of the sequences
we considered.

In two papers since then \cite{V7,V8}, we obtained results when $y_{0}=b^{2}$
for any positive integer $b$ and $-N_{\alpha}$ (defined below in Subsection~\ref{subsect:notation1})
is either a square or of the form $2^{\ell}p^{m}$ with $p$ an odd prime and $\ell$,
$m$ non-negative integers. In this paper, we extend these results to any $N_{\alpha}<0$.

\subsection{Notation}
\label{subsect:notation1}

We follow the same notation as in our previous papers on this subject. We recall
that notation here.

Let $a$, $b$ and $d$ be positive integers such that $d$ is not a square.
For $\alpha=a+b^{2}\sqrt{d}$, put $N_{\alpha}
=N_{\bbQ \left( \sqrt{d} \right)/\bbQ}(\alpha)=a^{2}-b^{4}d$ and let
$\varepsilon=\left( t+u\sqrt{d} \right)/2$ be a unit in
$\cO_{\bbQ \left( \sqrt{d} \right)}$ with $t$ and $u$ positive integers.

We define the two sequences $\left( x_{k} \right)_{k=-\infty}^{\infty}$
and $\left( y_{k} \right)_{k=-\infty}^{\infty}$ by
\begin{equation}
\label{eq:yk-defn}
x_{k}+y_{k}\sqrt{d}
=\alpha \varepsilon^{2k}.
\end{equation}

Observe that $x_{0}=a$,
$y_{0}=b^{2}$,
\begin{equation}
\label{eq:yPM1}
y_{1}= \frac{\left( b^{2}\left( t^{2}+du^{2} \right)+2atu \right)}{4}, \hspace*{3,0mm}
y_{-1}=\frac{\left( b^{2}\left( t^{2}+du^{2} \right)-2atu \right)}{4}
\end{equation}
and that both sequences satisfy
the recurrence relation
\begin{equation}
\label{eq:yk-recurrence}
u_{k+1}= \frac{t^{2}+du^{2}}{2} u_{k}-u_{k-1},
\end{equation}
for all $k \in \bbZ$.

To relate such sequences to the quadratic equations mentioned in the previous
subsection, observe that from \eqref{eq:yk-defn},
\[
x_{k}^{2}-dy_{k}^{2}=N_{\alpha}.
\]

We restrict the coefficient of $\sqrt{d}$ in $\alpha$ to being a square since
we are interested in squares among the $y_{k}$'s. We also choose $\alpha$ so
that $b^{2}$ is the smallest square among the $y_{k}$'s.

For any non-zero integer, $n$, let $\core(n)$ be the unique squarefree integer
such that $n/\core(n)$ is a square. We will put $\core(1)=1$.

\subsection{Conjectures}
\label{subsect:conjectures}

We first recall from \cite{V5} some conjectures about the number of distinct
squares in the sequence of $y_{k}$'s.

\begin{conjecture}
\label{conj:1-seq}
There are at most four distinct integer squares among the $y_{k}$'s.

If $\core \left( \left| N_{\alpha} \right| \right)|(2p)$ where $p$ is an odd
prime, then there
are at most three distinct integer squares among the $y_{k}$'s.

Furthermore, if $\left| N_{\alpha} \right|$ is a perfect square, then there are
at most two distinct integer squares among the $y_{k}$'s.
\end{conjecture}

A more general result than Conjecture~\ref{conj:1-seq} actually appears to be true.
Removing the restriction to even powers of $\varepsilon$ in \eqref{eq:yk-defn},
define $\left( x_{k}' \right)_{k=-\infty}^{\infty}$ and $\left( y_{k}' \right)_{k=-\infty}^{\infty}$
by
\begin{equation}
\label{eq:ykPrime-defn}
x_{k}'+y_{k}'\sqrt{d}
=\alpha \varepsilon^{k}.
\end{equation}

\begin{conjecture}
\label{conj:2-seq}
There are at most four distinct integer squares among the $y_{k}'$'s.

If $\left| N_{\alpha} \right|$ is a prime power or a perfect square,
then there are at most three distinct integer squares among the $y_{k}'$'s.
\end{conjecture}

Computational evidence for these conjectures was presented in Subsection~1.3 of \cite{V5}.

\subsection{Previous Results}
\label{subsect:prev-results}

The following is the main theorem in \cite{V5}. With the exceptions noted in parts~(a)
and (b), it establishes Conjecture~$\ref{conj:1-seq}$ when $b=1$ and $-N_{\alpha}$
is a square.

\begin{theorem}
\label{thm:b1-sqr}
Let $b=1$, $a$ and $d$ be positive integers, where $d$ is not a square, $N_{\alpha}<0$
and $-N_{\alpha}$ is a square.

\noindent
{\rm (a)} If $u=1$, $t^{2}-du^{2}=-4$, $N_{\alpha} \equiv 12 \pmod{16}$, $\gcd \left( a^{2},d \right)=1,4$
and one of $y_{\pm 1}$ is a perfect square, then there are at most three
distinct squares among the $y_{k}$'s.

\noindent
{\rm (b)} If $u=2$, $t^{2}-du^{2}=-4$, $N_{\alpha}$ is odd, $\gcd \left( a^{2},d \right)=1$
and one of $y_{\pm 1}$ is a perfect square, then there are at most three
distinct squares among the $y_{k}$'s.

\noindent
{\rm (c)} Otherwise, there are at most two distinct squares among the $y_{k}$'s.
\end{theorem}

The following is the main theorem in \cite{V6}.

\begin{theorem}
\label{thm:b1-pp}
Let $b=1$ and let $a$, $m$ and $p$ be non-negative integers with $a \geq 1$ and
$p$ a prime. Conjecture~$\ref{conj:1-seq}$ holds when $N_{\alpha}=-p^{m}$,
$-2p^{m}$, $-4p^{m}$, $-8p^{m}$ and $-16p^{m}$.
\end{theorem}

The following is Corollary~\corMainSqr{} in \cite{V7}. It extends Theorem~\ref{thm:b1-sqr}
to any $b$.

\begin{theorem}
\label{thm:any-b-sqr}
Let $a$, $b$ and $d$ be positive integers, where $d$ is not a square, $N_{\alpha}<0$
and $-N_{\alpha}$ is a square.

{\rm (a)} For $b=1,\ldots, 11$, there are at most five distinct squares among the
$y_{k}$'s.

{\rm (b)} For $b \geq 12$, there are at most five distinct squares among the $y_{k}$'s
provided that
\[
d \geq \frac{30 \left| N_{\alpha} \right|^{1/2}b^{26/11}}{u^{24/13}}.
\]
\end{theorem}

Lastly, the following is Corollary~\corMainPP{} in \cite{V8}. It extends Theorem~\ref{thm:b1-pp}
to any $b$.

\begin{theorem}
\label{thm:any-b-pp}
Let $a$, $b$ and $d$ be positive integers, where $d$ is not a square, $N_{\alpha}<0$
and $-N_{\alpha}=2^{\ell}p^{m}$ with $p$ an odd prime and $\ell$, $m$ non-negative
integers.

{\rm (a)} For $b=1,\ldots, 6$, there are at most seven distinct squares among the
$y_{k}$'s.

{\rm (b)} For $b \geq 7$, there are at most seven distinct squares among the $y_{k}$'s,
provided that
\[
d \geq \frac{65 \left| N_{\alpha} \right|^{1/2}b^{26/11}}{u^{24/13}}.
\]
\end{theorem}

\subsection{New Results}
\label{subsect:new-results}

In this paper, we make further progress towards establishing Conjecture~\ref{conj:1-seq}.
We obtain a small absolute upper bound for the number of distinct squares
in such sequences for $N_{\alpha}<0$, when $b$ is small or when $d$ is sufficiently
large.

Let $K$ be the largest negative integer such that $y_{K}>b^{2}$.

\begin{theorem}
\label{thm:1.2-seq-new}
Let $a$, $b$ and $d$ be positive integers, where $d$ is not a square and $N_{\alpha}<0$.
There are at most four distinct squares among the
$y_{k}$'s with $k \geq 3$ or $k \leq K-2$, and
\[
y_{k}> \frac{16 b^{4}\left| N_{\alpha} \right|^{4}}{\sqrt{d}}.
\]
\end{theorem}

\begin{theorem}
\label{thm:1.3-seq-new}
Let $a$, $b$ and $d$ be positive integers, where $d$ is not a square and $N_{\alpha}<0$.

{\rm (a)} For $b=1,2$ or $3$, there are at most four distinct squares among the
$y_{k}$'s with $k \geq 3$ or $k \leq K-2$.

{\rm (b)} For $b \geq 4$, there are at most four distinct squares among the
$y_{k}$'s with $k \geq 3$ or $k \leq K-2$, provided that
\[
d \geq
\frac{15\left| N_{\alpha} \right|^{3/4} b^{3/2}}{u^{3/2}}.
\]
\end{theorem}

Theorems~1.3 and 1.6 of \cite{V7} and \cite{V8} respectively are analogues of
this theorem for their respective values of $N_{\alpha}$. Theorems~\ref{thm:any-b-sqr}
and \ref{thm:any-b-pp} above are corollaries of these theorems in \cite{V7,V8}.

The following two corollaries follow immediately from Theorem~\ref{thm:1.3-seq-new}.

\begin{corollary}
\label{cor:1.3-seq-new}
Let $a$, $b$ and $d$ be positive integers, where $d$ is not a square, $N_{\alpha}<0$.

{\rm (a)} For $b=1,\ldots, 3$, there are at most $9$ distinct squares among the
$y_{k}$'s.

{\rm (b)} For $b \geq 4$, there are at most $9$ distinct squares among the $y_{k}$'s,
provided that
\[
d \geq \frac{15 \left| N_{\alpha} \right|^{3/4}b^{3/2}}{u^{3/2}}.
\]
\end{corollary}

\begin{corollary}
\label{cor:1.4-seq-new}
Let $a$, $b$ and $d$ be positive integers, where $d$ is not a square and $N_{\alpha}<0$.

If the smallest squares among the even-indexed
$y_{k}'$ and among the odd-indexed $y_{k}'$ are both at most $9$, then there are
at most $18$ distinct squares among the $y_{k}$'s.
\end{corollary}

The code used in this work is publicly available at \url{https://github.com/PV-314/hypgeom/any-b}.
The author is very happy to help interested readers who have any questions,
problems or suggestions for the use of this code.

\subsection{Future work}

We highlight here two areas where further work would lead to improvements in the results.

\vspace*{1.0mm}

\noindent
(1) Treatment of small squares among the $y_{k}$'s.
Many diophantine problems are solved by considering separately small and large
solutions. A technique to efficiently count small squares among the $y_{k}$'s
(i.e., squares among the $y_{k}$'s not satisfying the lower bound for $y_{k}$ in
Theorem~\ref{thm:1.2-seq-new}) would provide unconditional results like our
conjectures above.

\vspace*{1.0mm}

\noindent
(2) Improved treatment of the $r_{0}=1$ steps in our proof.
It is only in these two steps where we need to assume the existence of five sufficiently
large distinct squares, rather than three such squares. Moreover, in these steps,
we have simple explicit expressions for associated hypergeometric polynomials.
If these expressions, or other means, could be used to treat these steps better,
then we could prove an improved version of Theorem~\ref{thm:1.2-seq-new} with
only $7$ distinct squares.

\section{Lemmas about $\left( x_{k} \right)_{k=-\infty}^{\infty}$ and $\left( y_{k} \right)_{k=-\infty}^{\infty}$}
\label{sect:prelim}

The following lemma allows us to use the hypergeometric method.

\begin{lemma}
\label{lem:quad-rep}
Let $a \neq 0$, $b>0$ and $d$ be rational integers such that $d$ is not a square.
Put $\alpha=a+b^{2} \sqrt{d}$ and denote $N_{\bbQ(\sqrt{d})/\bbQ}(\alpha)$
by $N_{\alpha}$. Suppose that $N_{\alpha}$ is not a square, $x \neq 0$ and $y>0$
are rational integers with
\begin{equation}
\label{eq:quad-rep-assumption}
x+y^{2} \sqrt{d} = \alpha \epsilon^{2},
\end{equation}
where $\epsilon=\left( t+u\sqrt{d} \right)/2 \in \cO_{\bbQ \left( \sqrt{d} \right)}$
with $t$ and $u$ non-zero rational integers, norm $N_{\epsilon}=\pm 1$.

We can write
\begin{align}
f^{2} \left( x + N_{\epsilon}\sqrt{N_{\alpha}} \right)
&= \left( a+\sqrt{N_{\alpha}} \right) \left( r + s\sqrt{\core \left( N_{\alpha} \right)} \right)^{4}
\quad \text{and} \nonumber \\
\label{eq:fy-rel}
fy &= b \left( r^{2}-\core \left( N_{\alpha} \right)s^{2} \right),
\end{align}
for some integers $f$, $r$ and $s$ satisfying $f \neq 0$,
\[
f | \left( 4b^{2} \core \left( N_{\alpha} \right) \right) \text{ and }
0< f < 4b^{2} \sqrt{\core \left( N_{\alpha} \right)}.
\]
\end{lemma}

\begin{proof}
This is Proposition~\propRep{}(a) of \cite{V5} when $\epsilon$ is a unit in $\bbQ \left( \sqrt{d} \right)$.
\end{proof}

We will also need lower bounds for the elements in our sequences.

\begin{lemma}
\label{lem:Y-LB1}
Let the $y_{k}$'s be defined by $\eqref{eq:yk-defn}$ with the notation and assumptions
there. Suppose that $N_{\alpha}<0$.
Let $K$ be the largest negative integer such that $y_{K}>b^{2}$.

\noindent
{\rm (a)} Put $\overline{\alpha}=a-b^{2}\sqrt{d}$. We have
\begin{equation}
\label{eq:yLB1-gen}
y_{k}>
\left\{
\begin{array}{ll}
\dfrac{\alpha\varepsilon^{2k}}{2\sqrt{d}} & \text{for $k \geq 0$,} \\
\dfrac{-\overline{\alpha} \, \varepsilon^{2|k|}}{2\sqrt{d}} & \text{for $k<0$.}
\end{array}
\right.
\end{equation}

\noindent
{\rm (b)}
For all $k$, $2y_{k}$ is a positive integer.
The sequences $\left( y_{k} \right)_{k \geq 0}$ and 
$\left( y_{K+1}, y_{K}, y_{K-1}, y_{K-2}, \ldots \right)$ are increasing sequences of positive numbers.

\noindent
{\rm (c)} We have
\begin{equation}
\label{eq:yLB3-gen}
y_{k} \geq
\left\{
\begin{array}{ll}
\left( \left| N_{\alpha} \right|u^{2} / \left( 4b^{2} \right) \right) \left( 2du^{2}/5 \right)^{k-1} & \text{for $k>0$,} \\
\left( \left| N_{\alpha} \right|u^{2} / \left( 4b^{2} \right) \right) \left( 2du^{2}/5 \right)^{\max(0,K-k)} & \text{for $k<0$.}
\end{array}
\right.
\end{equation}
\end{lemma}

\begin{proof}
This is Lemma~\lemLB{} of \cite{V5}.
\end{proof}

Next we state a gap principle separating distinct squares in the sequence of
$y_{k}$'s.

For any $k$ with $y_{k}=y^{2}$, we let $r_{k}$ and $s_{k}$ be the
quantities in Lemma~\ref{lem:quad-rep}, using $\left( x_{k}, y_{k} \right)$
for $\left( x,y^{2} \right)$ there,
$\omega= \left( a-\sqrt{N_{\alpha}} \right)
/ \left( a+\sqrt{N_{\alpha}} \right)$ and letting $\zeta_{4}^{(k)}$ be the $4$-th
root of unity such that
\begin{equation}
\label{eq:zetaj-defn}
\left| \omega^{1/4} - \zeta_{4}^{(k)} \frac{r_{k}+s_{k}\sqrt{\core \left( N_{\alpha} \right)}}
{r_{k}-s_{k}\sqrt{\core \left( N_{\alpha} \right)}} \right|
= \min_{0 \leq \ell \leq 3} \left| \omega^{1/4} - e^{2\ell \pi i/4} \frac{r_{k}+s_{k}\sqrt{\core \left( N_{\alpha} \right)}}{r_{k}-s_{k}\sqrt{\core \left( N_{\alpha} \right)}} \right|.
\end{equation}

\begin{lemma}
\label{lem:gap}
Let the $y_{k}$'s be defined as in $\eqref{eq:yk-defn}$ with $N_{\alpha}<0$.

If $y_{k_{1}}$, $y_{k_{2}}$ and $y_{k_{3}}$
are three distinct squares with none of $k_{1}$, $k_{2}$, $k_{3}$ between $K+1$ and $0$,
inclusive, and
$y_{k_{3}}>y_{k_{2}}>y_{k_{1}} \geq \max \left\{ 4\sqrt{\left| N_{\alpha} \right|/d},
\left( 16b^{2} \left| N_{\alpha} \right|^{2} / d \right)^{2}/60 \right\}$,
then there exist distinct $i,j \in \{ k_{1},k_{2}, k_{3} \}$ such that
\[
y_{j} > 1.43 \frac{d}{b^{2} \left| N_{\alpha} \right|^{2}} y_{i}^{5/2}.
\]

Moreover, we can let $i$ and $j$ be any distinct elements of $\{ k_{1},k_{2}, k_{3} \}$
such that $\zeta_{4}^{(i)}= \pm \zeta_{4}^{(j)}$, where $\zeta_{4}^{(i)}$ and
$\zeta_{4}^{(j)}$ are as defined in $\eqref{eq:zetaj-defn}$.
\end{lemma}

\begin{remark-nonum}
The purpose of the $1/60$ in the second term of the max is to reduce the size
of the constant in \eqref{eq:dUB-from-ass-2b-b} in Subsection~\ref{subsect:proof-thm13b}.
This will reduce the amount of computation that is necessary in Subsection~\ref{subsect:proof-thm13a}
for the proof of Theorem~\ref{thm:1.3-seq-new}(a).
\end{remark-nonum}

\begin{proof}
We will apply Lemma~\lemGap{}(b) of \cite{V5} and its proof. Although it is assumed
in Lemma~\lemGap{}(b) of \cite{V5} that $-N_{\alpha}$ is not a square, this assumption
is not used in the proof and it is only added in the statement of part~(b) to distinguish
it from part~(a).

From Lemma~\lemGap{}(b) of \cite{V5}, there exist distinct
$i,j \in \{ k_{1},k_{2}, k_{3} \}$ such that
\[
y_{j} > 15.36 \left( \frac{b^{2}d}{f_{i}f_{j}\left| N_{\alpha} \right|} \right)^{2}
y_{i}^{3}.
\]

By Lemma~\ref{lem:quad-rep}, we have
\begin{equation}
\label{eq:fm1-fm3-ub}
f_{i} f_{j} \leq 16b^{4} \left| N_{\alpha} \right|,
\end{equation}
so
\[
y_{j}
> 15.36 \left( \frac{b^{2}d}{f_{i}f_{j}\left| N_{\alpha} \right|} \right)^{2} y_{i}^{3}
\geq 15.36 \left( \frac{d}{16b^{2} \left| N_{\alpha} \right|^{2}} \right)^{2} y_{i}^{3}.
\]

Applying the second value in the max in the lower bound for $y_{k_{1}}$ in the
statement of this lemma, we obtain
\[
y_{j} > 0.256y_{i}^{2}.
\]

The condtion that none of $k_{1}$, $k_{2}$ or $k_{3}$ are between $K+1$ and
$0$ implies that $y_{i}$ is a square greater than $b^{2}$. So $y_{i} \geq (b+1) \geq 4$.
Thus
\[
y_{j} > 0.256y_{i}^{2} \geq 1.024y_{i}.
\]

We can now proceed as in the proof of Lemma~\lemGap{}(b) of \cite{V5} from
equation~(\eqGPb{}) onwards, which we restate here:
\[
\frac{2b}{\sqrt{f_{i}f_{j}} \left( y_{i}y_{j} \right)^{1/4}}
< 0.5051b^{2} \sqrt{\frac{\left| N_{\alpha} \right|}{a^{2}+\left| N_{\alpha} \right|}}
\left( \frac{1}{y_{i}} + \frac{1}{y_{j}} \right).
\]

Applying $y_{j}>1.024y_{i}$, we obtain
\[
\frac{2b}{\sqrt{f_{i}f_{j}} \left( y_{i}y_{j} \right)^{1/4}}
< 0.99837\sqrt{\frac{\left| N_{\alpha} \right|}{d}}
\frac{1}{y_{i}}.
\]

This gives
\[
y_{j} > 16.1\left( \frac{b^{2}d}{f_{i}f_{j} \left| N_{\alpha} \right|} \right)^{2} y_{i}^{3}.
\]

Using Maple, we repeated this process $14$ more times, yielding
\[
y_{j} > 177.77\left( \frac{b^{2}d}{f_{i}f_{j} \left| N_{\alpha} \right|} \right)^{2} y_{i}^{3}.
\]

Applying the second value in the max in the lower bound for $y_{k_{1}}$
and \eqref{eq:fm1-fm3-ub}, we obtain
\[
y_{j} > \frac{177.77}{\sqrt{60}} \frac{b^{2}d}{f_{i}f_{j} \left| N_{\alpha} \right|} y_{i}^{5/2}.
\]

Applying \eqref{eq:fm1-fm3-ub} once more:
\[
y_{j} > \frac{177.77}{\sqrt{60}} \frac{b^{2}d}{16b^{4} \left| N_{\alpha} \right|^{2}} y_{i}^{5/2}
 > 1.43 \frac{d}{b^{2} \left| N_{\alpha} \right|^{2}} y_{i}^{5/2}.
\]

The statement about
$\zeta_{4}^{(i)}= \pm \zeta_{4}^{(j)}$ is proven in the second paragraph of the
proof of Lemma~\lemGap{}(b) of \cite{V5} (see page~317 in \cite{V5}).
\end{proof}

\section{Proof of Theorem~\ref{thm:1.2-seq-new}}
\label{sect:proof-thm12}

\subsection{Choice of indices}
\label{subsect:k-and-ell}

To prove Theorem~\ref{thm:1.2-seq-new}, we assume there are five distinct squares,
\[
y_{k_{5}}>y_{k_{4}}>y_{k_{3}}>y_{k_{2}}>y_{k_{1}},
\]
satisfying
\begin{equation}
\label{eq:5-sqrs}
k_{1} \geq 3 \text{ or } k_{1} \leq K-2,
\end{equation}
recalling that $K$ is the largest negative integer such that $y_{K} > b^{2}$.

We also need to ensure that the conditions in Lemmas~\ref{lem:gap} above and
Lemma~\lemTheta{} of \cite{V5} hold. This will aid us in our choice of indices below, as well as in our application
of our gap principle in Subsections~\ref{subsect:thm12-step-i} and \ref{subsect:thm12-step-ii}.
So we will assume that
\begin{equation}
\label{eq:yks-LB1}
y_{k_{5}}>y_{k_{4}}>y_{k_{3}}>y_{k_{2}}>y_{k_{1}}
\geq \max \left\{ 4\sqrt{\left| N_{\alpha} \right|/d}, (64/15)b^{4} \left| N_{\alpha} \right|^{4}/d^{2} \right\}
\end{equation}

Furthermore, to obtain our lower bounds for $E$ and $Q$ in \eqref{eq:E-LB1} and
\eqref{eq:Q-LB1}, we will also require that
\begin{equation}
\label{eq:yk1-LB1}
y_{k_{1}} \geq \max \left\{ 4\left| N_{\alpha} \right|/\sqrt{d}, |g|\cN_{d',4}/\sqrt{d} \right\}.
\end{equation}

In Subsection~\subsectPrereqs{} of \cite{V7}, we showed that \eqref{eq:yk1-LB1}
holds if $y_{k_{1}} \geq 4\left| N_{\alpha} \right|/\sqrt{d}$.
Combining this with \eqref{eq:yks-LB1}, we find that \eqref{eq:yks-LB1} and
\eqref{eq:yk1-LB1} both hold if
\begin{equation}
\label{eq:yks-LB2}
y_{k_{5}}>y_{k_{4}}>y_{k_{3}}>y_{k_{2}}>y_{k_{1}}
\geq \max \left\{ 4\left| N_{\alpha} \right|/\sqrt{d}, (64/15)b^{4}\left| N_{\alpha} \right|^{4}/d^{2} \right\}.
\end{equation}

Using \eqref{eq:zetaj-defn}, we have $\zeta_{4}^{(k_{1})}$, \ldots, $\zeta_{4}^{(k_{5})}$
associated with $k_{1}$,\ldots, $k_{5}$. At least three of these must be either
all equal to $\pm 1$ or else all equal to $\pm i$. We label
the three associated indices as $m_{1}$, $m_{2}$ and $m_{3}$, where
$m_{1}<m_{2}<m_{3}$. If the set of at least three elements contains more than
three elements, we can choose any three among them.

It is in order to have three such elements that we assume we have five squares
in \eqref{eq:yks-LB1}. The reason we want three such elements is because in
Subsections~\ref{subsect:thm12-step-i} and \ref{subsect:thm12-step-ii}, we will
apply the gap principle in Lemma~\ref{lem:gap} twice.

We put
$\omega_{m_{1}}=\left( x_{m_{1}}+N_{\varepsilon^{m_{1}}}\sqrt{N_{\alpha}} \right)/\left( x_{m_{1}}-N_{\varepsilon^{m_{1}}}\sqrt{N_{\alpha}} \right)$
and let $\zeta_{4}$ be the $4$-th root of unity such that
\[
\left| \omega_{m_{1}}^{1/4} - \zeta_{4} \frac{x-y\sqrt{\core \left( N_{\alpha} \right)}}{x+y\sqrt{\core \left( N_{\alpha} \right)}} \right|
= \min_{0 \leq j \leq 3} \left| \omega_{m_{1}}^{1/4} - e^{2j\pi i/4} \frac{x-y\sqrt{\core \left( N_{\alpha} \right)}}{x+y\sqrt{\core \left( N_{\alpha} \right)}} \right|,
\]
where $x+y\sqrt{\core \left( N_{\alpha} \right)}=\left(  r_{m_{1}}-s_{m_{1}}\sqrt{\core \left( N_{\alpha} \right)} \right)
\left( r_{m_{3}}+s_{m_{3}}\sqrt{\core \left( N_{\alpha} \right)} \right)$ with
$\left( r_{m_{1}}, s_{m_{1}} \right)$ and $\left( r_{m_{3}}, s_{m_{3}} \right)$
as in Lemma~\ref{lem:quad-rep},
which are associated with $\left( x_{m_{1}}, y_{m_{1}} \right)$ and
$\left( x_{m_{3}}, y_{m_{3}} \right)$, respectively.
We can take $\zeta_{4} \in \bbQ \left( \sqrt{\core \left( N_{\alpha} \right)} \right)$.
This is immediate since $N_{\alpha}$ is a perfect square. Otherwise, it follows
by applying Lemma~\lemTheta{}(b) in \cite{V5} with $k=m_{1}$ and $\ell=m_{3}$, along with
the fact that we chose $m_{1}$, $m_{2}$ and $m_{3}$ so that
$\zeta^{(m_{1})}=\pm \zeta^{(m_{3})}$, we find that $\zeta_{4}=\pm 1 \in \bbQ \left( \sqrt{\core \left( N_{\alpha} \right)} \right)$.

This is important for us here as $\zeta_{4} \left( x-y\sqrt{\core \left( N_{\alpha} \right)} \right)
/ \left( x+y\sqrt{\core \left( N_{\alpha} \right)} \right)$ must be in an
imaginary quadratic field in order to apply Lemma~\lemDio{} in \cite{V5}
to obtain a lower bound for the rightmost quantity in \eqref{eq:29} below.

\subsection{Further Prerequisites}
\label{subsect:preq}

In this subsection, we collect some inequalities what will be required in the
subsections that follow.

Recall equation~(\eqOmegaApproxUB{}) from \cite{V5} with $\ell=m_{3}$ and our
introduction here of a subscript on $\omega$:
\begin{equation}
\label{eq:29}
\frac{2\sqrt{\left| N_{\alpha} \right|}}{\sqrt{d} \, y_{m_{3}}}
= \left| \omega_{m_{1}} - \left( \frac{x-y\sqrt{\core \left( N_{\alpha} \right)}}{x+y\sqrt{\core \left( N_{\alpha} \right)}} \right)^{4} \right|
> 3.959 \left| \omega_{m_{1}}^{1/4} - \zeta_{4} \frac{x-y\sqrt{\core \left( N_{\alpha} \right)}}{x+y\sqrt{\core \left( N_{\alpha} \right)}} \right|.
\end{equation}

We need to derive a lower bound for the far-right quantity in \eqref{eq:29}.
To do so, we shall use the lower bounds in Lemma~\lemDio{} in \cite{V5} with
a sequence of good approximations $p_{r}/q_{r}$ obtained from the hypergeometric
functions. So we collect here the required quantities.

Since $y_{m_{1}} \geq 4\sqrt{\left| N_{\alpha} \right|/d}$ (from \eqref{eq:yks-LB2}),
we obtain
\begin{equation}
\label{eq:25}
x_{m_{1}}^{2} = dy_{m_{1}}^{2}+N_{\alpha}
=dy_{m_{1}}^{2} \left( 1 + \frac{N_{\alpha}}{dy_{m_{1}}^{2}} \right)
\geq 0.9375dy_{m_{1}}^{2}.
\end{equation}

So
\begin{equation}
\label{eq:26}
\sqrt{x_{m_{1}}^{2}-N_{\alpha}}= \sqrt{dy_{m_{1}}^{2}}
< 1.04 x_{m_{1}}.
\end{equation}

Using the notation from Subsection~\subSectConstr{} of \cite{V5}, let
$t'=\core \left( N_{\alpha} \right)$,
$u_{1}=2x_{m_{1}}$, $u_{2}=2\sqrt{N_{\alpha}/\core \left( N_{\alpha} \right)}$ and
$d'$ is as defined in equation~(\eqDPrimeDef{}) of \cite{V5}.

Recall equation~(\eqELB{}) in \cite{V5} with $k$ there being $m_{1}$ here:
\begin{equation}
\label{eq:E-LB1}
E > \frac{0.1832|g|\cN_{d',4}\sqrt{d} \, y_{m_{1}}}{\left| N_{\alpha} \right|}.
\end{equation}

From Lemma~\ref{lem:Y-LB1}(c), $|g|\cN_{d',4} \geq 2$ and $y_{m_{1}} \geq 4\left| N_{\alpha} \right|/\sqrt{d}$,
we have $E>1$, as required for its use with Lemma~\lemDio{} in \cite{V5}.

Similarly, using \eqref{eq:25} and Lemma~\ref{lem:Y-LB1}(c), we have
\begin{equation}
\label{eq:Q-LB1}
Q > \frac{2e^{1.68}\left( 1+\sqrt{0.9375} \right) \sqrt{d}y_{m_{1}}}{|g|\cN_{d',4}}
>\frac{21.12\sqrt{d} \, y_{m_{1}}}{|g|\cN_{d',4}}
\geq 21.12,
\end{equation}
where we use $y_{m_{1}} \geq |g|\cN_{d',4}/\sqrt{d}$
to establish the last inequality.
For \eqref{eq:E-LB1}\footnote{PV TODO: identify where} and \eqref{eq:Q-LB1}, we
have used \eqref{eq:yk1-LB1} (and hence \eqref{eq:yks-LB2}).

From equation~(\eqQUB{}) in \cite{V5}, we have
\begin{equation}
\label{eq:Q-UB2}
Q <\frac{21.47\sqrt{d} \, y_{m_{1}}}{|g|\cN_{d',4}}.
\end{equation}

Writing $\omega_{m_{1}}=e^{i\varphi_{m_{1}}}$, with $-\pi < \varphi_{m_{1}} \leq \pi$,
from equation~(\eqEllUB{}) in \cite{V5}, we can take
\begin{equation}
\label{eq:ell-UB}
\ell_{0}<0.458\sqrt{\left| N_{\alpha} \right|}/\left| x_{m_{1}} \right|.
\end{equation}

Also from Lemma~\lemTheta{}(a) in \cite{V5}, we have $\left| \varphi_{m_{1}} \right|<0.6$,
so the condition $|\omega-1|<1$ required in Lemma~\lemHypg{}\, of \cite{V5} to
apply the hypergeometric method is satisfied.

Let $q=x+y\sqrt{\core \left( N_{\alpha} \right)}=\left( r_{m_{1}}-s_{m_{1}}\sqrt{\core \left( N_{\alpha} \right)} \right)
\left( r_{m_{3}}+s_{m_{3}}\sqrt{\core \left( N_{\alpha} \right)} \right)$
and $p=x-y\sqrt{\core \left( N_{\alpha} \right)}$. Recall from (\eqQAbs{}) in
\cite{V5} that
\begin{equation}
\label{eq:qAbs}
|q| = \frac{\sqrt{f_{m_{1}}f_{m_{3}}} \left( y_{m_{1}}y_{m_{3}} \right)^{1/4}}{b}.
\end{equation}

We know from \eqref{eq:5-sqrs} that $m_{1} \geq 3$ or $m_{1} \leq K-2$, so from
Lemma~\ref{lem:Y-LB1}(c) we obtain
\begin{equation}
\label{eq:ym1-LB}
y_{m_{1}} \geq \frac{|N_{\alpha}|d^{2}u^{6}}{25b^{2}}.
\end{equation}

We are now ready to deduce the required contradiction from the assumptions above.
We define $r_{0}$ as in Lemma~\lemDio{} of \cite{V5} to be the smallest positive
integer such that
$\left( Q-1/E  \right)\ell_{0}|B|/\left( Q-1 \right)<cE^{r_{0}}$. We will
take $c=0.75$ here, as in \cite{V5}.
In each of the four main steps below, we will obtain an upper bound for $y_{m_{1}}$
required for that step to hold. In this way, we will show that none of the four
main steps can hold for $y_{m_{1}}$ sufficiently large.

\subsection{$r_{0}=1$ and $\zeta_{4}p/q \neq p_{1}/q_{1}$ for all $4$-th roots of unity, $\zeta_{4}$}
\label{subsect:thm12-step-i}

We start with the upper bound for $y_{m_{3}}$ in equation~(\eqYUBStepOne{}) in
\cite{V7}:
\[
y_{m_{3}}^{3} < 545(b(1-c))^{-4} \left( N_{\alpha}f_{m_{1}}f_{m_{3}} \right)^{2}y_{m_{1}}^{5}.
\]

Applying \eqref{eq:fm1-fm3-ub}, we obtain
\begin{equation}
\label{eq:y3UB}
y_{m_{3}}^{3} < 140000(1-c)^{-4} N_{\alpha}^{4}b^{4}y_{m_{1}}^{5}.
\end{equation}

Applying Lemma~\ref{lem:gap} twice, and \eqref{eq:fm1-fm3-ub} for any non-zero
$i$ and $j$, we have
\[
y_{m_{3}}
> \left( \frac{1.43d}{b^{2} \left| N_{\alpha} \right|^{2}} \right) y_{m_{2}}^{5/2}
> \left( \frac{1.43d}{b^{2} \left| N_{\alpha} \right|^{2}} \right)^{7/2} y_{m_{1}}^{25/4}.
\]
This is possible since, by our choice of $m_{1}$, $m_{2}$ and $m_{3}$,
we have $\zeta^{(m_{1})}$, $\zeta^{(m_{2})}$ and $\zeta^{(m_{3})}$
all equal to $\pm 1$ times each other.

Combining the lower bound that this provides for $y_{m_{3}}^{3}$ with the upper
bound for $y_{m_{3}}^{3}$ in \eqref{eq:y3UB} and cancelling the common factor
of $y_{m_{1}}^{5}$ on both sides, we find that
\[
\left( \frac{1.43d}{b^{2} \left| N_{\alpha} \right|^{2}} \right)^{21/2}y_{m_{1}}^{55/4}
< 140000(1-c)^{-4}N_{\alpha}^{4}b^{4},
\]
provided that \eqref{eq:yks-LB1} holds, which holds if \eqref{eq:yks-LB2} holds.

We need \eqref{eq:5-sqrs} here, as it implies that $m_{1}, m_{3} \neq 0$, which
is not permitted in Lemma~\ref{lem:gap} (which we use above).

Rearranging and using $c=0.75$ (a choice we justified in Subsection~\subSectStepiii{} of \cite{V5}),
we see that
\begin{equation}
\label{eq:yk-UB-step-i}
\frac{2.7 b^{20/11}\left| N_{\alpha} \right|^{20/11}}{d^{42/55}}
> y_{m_{1}}
\end{equation}
must hold here.

\subsection{$r_{0}=1$ and $\zeta_{4}p/q = p_{1}/q_{1}$ for some $4$-th root of unity, $\zeta_{4}$}
\label{subsect:thm12-step-ii}

As in Subsection~\ref{subsect:thm12-step-i}, we start by calling an upper bound for
$y_{m_{3}}$ that holds for all $r_{0} \geq 1$ with $\zeta_{4}p/q=p_{r_{0}}/q_{r_{0}}$
for some $4$-th root of unity, $\zeta_{4}$. From (\eqYLUBUsingYK{}) in \cite{V5}
with $k$ and $\ell$ there set to $m_{1}$ and $m_{3}$, respectively, we have
\begin{equation}
\label{eq:yEll-UB-step2}
1.73r_{0}^{1/2} \left(4\frac{d}{\left| N_{\alpha} \right|} \right)^{r_{0}}y_{m_{1}}^{2r_{0}+1}
>y_{m_{3}}.
\end{equation}

Specialising this inequality to the case of $r_{0}=1$ and applying the gap principle
in Lemma~\ref{lem:gap} twice, we obtain
\[
6.92 \frac{d}{\left| N_{\alpha} \right|} y_{m_{1}}^{3}>y_{m_{3}}
> \left( \frac{1.43d}{b^{2} \left| N_{\alpha} \right|^{2}} \right) y_{m_{2}}^{5/2}
> \left( \frac{1.43d}{b^{2} \left| N_{\alpha} \right|^{2}} \right)^{7/2} y_{m_{1}}^{25/4},
\]
provided that
$y_{m_{1}} \geq \max \left\{ 4\sqrt{\left| N_{\alpha} \right|/d}, (64/15)b^{4} \left| N_{\alpha} \right|^{4} / d^{2} \right\}$
(which holds by assumption \eqref{eq:yks-LB2} above).

Rearranging this, we find that
\begin{equation}
\label{eq:yk-UB-step-ii}
\frac{1.24 \left| N_{\alpha} \right|^{24/13}b^{28/13}}{d^{10/13}}>y_{m_{1}}
\end{equation}
must hold.

\subsection{$r_{0}>1$, $\zeta_{4}p/q \neq p_{r_{0}}/q_{r_{0}}$ for all $4$-th
roots of unity, $\zeta_{4}$}
\label{subsect:thm12-step-iii}

Recall equation~(\eqYUBStepThree{}) from \cite{V7}:
\begin{equation}
\label{eq:y1UB-step3-b}
\left( f_{m_{1}}f_{m_{3}} \right)^{8}
\frac{0.0047\left| N_{\alpha} \right|^{2} \left( |g|\cN_{d',4} \right)^{4}}{(bd)^{4}}
> \left( \frac{y_{m_{1}}^{4}d^{2}\left( |g|\cN_{d',4} \right)^{8}}
{59.2^{4} \left| N_{\alpha} \right|^{6}} \right)^{2r_{0}-1}.
\end{equation}

Applying \eqref{eq:fm1-fm3-ub}, the left-hand side of \eqref{eq:y1UB-step3-b}
is at most
\begin{equation}
\label{eq:yk-UB-step-iiia}
b^{28}\frac{2.02 \cdot 10^{7} \left| N_{\alpha} \right|^{10} \left( |g|\cN_{d',4} \right)^{4}}{d^{4}}.
\end{equation}

If
\begin{equation}
\label{eq:yk-UB-step-iiib}
\frac{59.2 \cdot 67.1^{1/(2r_{0}-1)}\left| N_{\alpha} \right|^{3/2+5/(2(2r_{0}-1))} b^{7/(2r_{0}-1)}}
{d^{1/2+1/(2r_{0}-1)} \left( |g|\cN_{d',4} \right)^{2-1/(2r_{0}-1)}}
>y_{m_{1}},
\end{equation}
then \eqref{eq:yk-UB-step-iiia} is greater than the right-hand side of \eqref{eq:y1UB-step3-b}.
We now simplify this bound. First, from Lemma~\lemGNd{} in \cite{V5}, we know that
$|g|\cN_{d',4} \geq 2$. From this and $r_{0} \geq 2$, we have
$\left( |g|\cN_{d',4} \right)^{-2+1/(2r_{0}-1)} \leq
\left( |g|\cN_{d',4} \right)^{1/(2r_{0}-1)}/4$. Applying this, $d \geq 2$ and
$r_{0} \geq 2$, we obtain
\begin{align*}
& \frac{59.2 \cdot 67.1^{1/(2r_{0}-1)}\left| N_{\alpha} \right|^{3/2+5/(2(2r_{0}-1))} b^{7/(2r_{0}-1)}}
{d^{1/2+1/(2r_{0}-1)} \left( |g|\cN_{d',4} \right)^{2-1/(2r_{0}-1)}}
\leq
\frac{59.2 b^{7/(2r_{0}-1)} \left| N_{\alpha} \right|^{3/2+5/(2(2r_{0}-1))}}
{4d^{1/2}} \cdot \left( \frac{67.1 \cdot 2}{d} \right)^{1/(2r_{0}-1)} \\
<& \frac{61b^{7/3} \left| N_{\alpha} \right|^{7/3}}{d^{1/2}}.
\end{align*}

So
\begin{equation}
\label{eq:yk-UB-step-iii}
\frac{61b^{7/3} \left| N_{\alpha} \right|^{7/3}}{d^{1/2}}>y_{m_{1}}
\end{equation}
must hold.

\subsection{$r_{0}>1$ and $\zeta_{4}p/q = p_{r_{0}}/q_{r_{0}}$ for some $4$-th root of unity, $\zeta_{4}$}
\label{subsect:thm12-step-iv}

Recall equation~(\eqUBStepFour{}) from \cite{V7}:
\[
1
> \frac{0.607b^{4}d}{f_{m_{1}}^{2}f_{m_{3}}^{2} \left| N_{\alpha} \right|}
  \left( \frac{0.0002344 |g|^{4}\cN_{d',4}^{4}}{\left| N_{\alpha} \right|^{3}} dy_{m_{1}}^{2} \right)^{r_{0}-1}.
\]

We now proceed similarly to the way we did in Subsection~\ref{subsect:thm12-step-iii}.

We have $f_{m_{1}}f_{m_{3}} \leq 16b^{4} \left| N_{\alpha} \right|$ from
\eqref{eq:fm1-fm3-ub}, so
\[
\frac{422b^{4}\left| N_{\alpha} \right|^{3}}{d}
> \left( \frac{0.0002344 |g|^{4}\cN_{d',4}^{4}}{\left| N_{\alpha} \right|^{3}} dy_{m_{1}}^{2} \right)^{r_{0}-1}.
\]

Rearranging this inequality, we find that it holds if
\begin{equation}
\label{eq:yk-UB-step-iva}
65.32\frac{b^{2/(r_{0}-1)}\left| N_{\alpha} \right|^{3/2+3/(2(r_{0}-1))} \cdot 422^{1/(2(r_{0}-1))}}
{|g|^{2}\cN_{d',4}^{2} d^{1/2+1/(2(r_{0}-1))}}>y_{m_{1}}.
\end{equation}

Applying $|g|\cN_{d',4} \geq 2$ (from Lemma~\lemGNd{} in \cite{V5}),
$d \geq 2$ and $r_{0} \geq 2$, we obtain
\begin{align*}
& 65.32\frac{b^{2/(r_{0}-1)}\left| N_{\alpha} \right|^{3/2+3/(2(r_{0}-1))} \cdot 422^{1/(2(r_{0}-1))}}
{|g|^{2}\cN_{d',4}^{2} d^{1/2+1/(2(r_{0}-1))}}
\leq 16.33\frac{b^{2/(r_{0}-1)}\left| N_{\alpha} \right|^{3/2+3/(2(r_{0}-1))}}
{d^{1/2}} \left( \frac{422}{d} \right)^{1/(2(r_{0}-1))}\\
<& \frac{238 b^{2}\left| N_{\alpha} \right|^{3}}{d^{1/2}}.
\end{align*}

Hence the above upper bound for $y_{m_{1}}$ holds, if
\begin{equation}
\label{eq:yk-UB-step-iv}
\frac{238 b^{2}\left| N_{\alpha} \right|^{3}}{d^{1/2}}>y_{m_{1}}.
\end{equation}

\subsection{Consolidation}
\label{subsect:thm12-consolidation}

Recall that we have assumed that there are five distinct squares,
$y_{k_{5}}>y_{k_{4}}>y_{k_{3}}>y_{k_{2}}>y_{k_{1}}$ with $k_{1} \geq 3$ or
$k_{1} \leq K-2$. Also, by \eqref{eq:yks-LB2}, we have
\[
y_{k_{1}} \geq \max \left\{ 4\left| N_{\alpha} \right|/\sqrt{d}, (64/15)b^{4}\left| N_{\alpha} \right|^{4}/d^{2} \right\}.
\]
But from \eqref{eq:yk-UB-step-i}, \eqref{eq:yk-UB-step-ii}, \eqref{eq:yk-UB-step-iii}
and \eqref{eq:yk-UB-step-iv}, we see that if we also have
\[
y_{m_{1}} \geq y_{k_{1}} > \max \left\{
\begin{array}{l}
\displaystyle\frac{2.7 b^{20/11}\left| N_{\alpha} \right|^{20/11}}{d^{42/55}},
\displaystyle\frac{1.3 \left| N_{\alpha} \right|^{24/13}b^{28/13}}{d^{10/13}},
\displaystyle\frac{61b^{7/3} \left| N_{\alpha} \right|^{7/3}}{d^{1/2}},
\displaystyle\frac{238 b^{2}\left| N_{\alpha} \right|^{3}}{d^{1/2}}
\end{array}
\right\},
\]
then we cannot be in any
of the four cases in Subsections~\ref{subsect:thm12-step-i}--\ref{subsect:thm12-step-iv}.
Hence we get a contradiction to our assumption that there are five such distinct
squares.

In case a more precise version of Theorem~\ref{thm:1.2-seq-new} is helpful, we
state this result here as a proposition.

\begin{proposition}
\label{prop:1.2-seq-new}
Let $a$, $b$ and $d$ be positive integers, where $d$ is not a square and
$N_{\alpha}<0$. There are at most four distinct squares among the
$y_{k}$'s with $k \geq 3$ or $k \leq K-2$, and
\[
y_{k} > \max \left\{
\frac{2.7 b^{20/11}\left| N_{\alpha} \right|^{20/11}}{d^{42/55}},
\frac{1.3 \left| N_{\alpha} \right|^{24/13}b^{28/13}}{d^{10/13}},
\frac{61b^{7/3} \left| N_{\alpha} \right|^{7/3}}{d^{1/2}},
\frac{238 b^{2}\left| N_{\alpha} \right|^{3}}{d^{1/2}},
\frac{4\left| N_{\alpha} \right|}{d^{1/2}},
\frac{64b^{4}\left| N_{\alpha} \right|^{4}}{15d^{2}}
\right\}.
\]
\end{proposition}

Taking the largest exponent on $\left| N_{\alpha} \right|$, the largest exponent
on $b$ and the smallest exponent on $d$, we see that this lower bound for $y_{k}$
is satisfied if $y_{k}>cb^{4}\left| N_{\alpha} \right|^{4}/\sqrt{d}$ for some
$c>0$. We may also assume that $N_{\alpha} \leq -15$ ($-14 \leq N_{\alpha} \leq -1$
are covered by \cite{V7,V8}). From this and $d \geq 2$, we find that the largest
constant obtained from the above lower bound for $y_{k}$ is $238/15^{4-3}=15.86\ldots$.
Hence if
\[
y_{k_{1}}>\frac{16 b^{4}\left| N_{\alpha} \right|^{4}}{\sqrt{d}},
\]
then the above lower bounds for $y_{k_{1}}$ hold, proving Theorem~\ref{thm:1.2-seq-new}.

\section{Proof of Theorem~\ref{thm:1.3-seq-new}}
\label{sect:proof-thm13}

Here we will combine our upper bounds for $y_{m_{1}}$ in Section~\ref{sect:proof-thm12}
with \eqref{eq:ym1-LB} to obtain Theorem~\ref{thm:1.3-seq-new}. As in
Section~\ref{sect:proof-thm12}, we proceed in separate steps.

Recall throughout this section that we have assumed there are three squares
among the $y_{k}$'s with $k \geq 3$ or $k \leq K-2$.

\subsection{$r_{0}=1$ and $\zeta_{4}p/q \neq p_{1}/q_{1}$ for all $4$-th roots of unity, $\zeta_{4}$}
\label{subsect:thm13-step-i}

In Subsection~\ref{subsect:thm12-step-i} (see \eqref{eq:yk-UB-step-i}), we showed
that
\[
y_{m_{1}} < \frac{2.7 b^{20/11}\left| N_{\alpha} \right|^{20/11}}{d^{42/55}}.
\]

This is not possible if
\[
\frac{2.7 b^{20/11}\left| N_{\alpha} \right|^{20/11}}{d^{42/55}}
<\frac{|N_{\alpha}|d^{2}u^{6}}{25b^{2}}.
\]
This holds if
\begin{equation}
\label{eq:thm13-step-i-dLB1}
d>\frac{4.6 \left| N_{\alpha} \right|^{45/152} b^{105/76}}{u^{165/76}}.
\end{equation}
It is useful to get a lower bound for $d$ just in terms of $b$ and $u$.
Applying $\left| N_{\alpha} \right| < db^{4}$, we also see that the upper
and lower bounds for $y_{m_{1}}$ in \eqref{eq:yk-UB-step-i} and \eqref{eq:ym1-LB}
contradict each other if
\begin{equation}
\label{eq:thm13-step-i-dLB2}
d>\frac{8.8b^{390/107}}{u^{330/107}}.
\end{equation}

Hence we cannot be in this step if either \eqref{eq:thm13-step-i-dLB1} or \eqref{eq:thm13-step-i-dLB2}
holds, along with the earlier assumptions made in Subsection~\ref{subsect:k-and-ell}.

\subsection{$r_{0}=1$ and $\zeta_{4}p/q=p_{1}/q_{1}$ for some $4$-th root of unity, $\zeta_{4}$}
\label{subsect:thm13-step-ii}

Combining the upper bound for $y_{m_{1}}$ in \eqref{eq:yk-UB-step-ii} with the
lower bound for it in \eqref{eq:ym1-LB}, we get a contradiction if
\[
\frac{1.24 \left| N_{\alpha} \right|^{24/13}b^{28/13}}{d^{10/13}}
< \frac{|N_{\alpha}|d^{2}u^{6}}{25b^{2}}.
\]

This inequality holds if
\begin{equation}
\label{eq:thm13-step-ii-dLB1}
d > \frac{3.46\left| N_{\alpha} \right|^{11/36} b^{3/2}}{u^{13/6}}.
\end{equation}

Applying $\left| N_{\alpha} \right|<db^{4}$, we obtain
\begin{equation}
\label{eq:thm13-step-ii-dLB2}
d > \frac{6b^{98/25}}{u^{78/25}}.
\end{equation}

\subsection{$r_{0}>1$, $\zeta_{4}p/q \neq p_{r_{0}}/q_{r_{0}}$ for all $4$-th
roots of unity, $\zeta_{4}$}
\label{subsect:thm13-step-iii}

Combining \eqref{eq:yk-UB-step-iiib} with \eqref{eq:ym1-LB},
\[
\frac{|N_{\alpha}|d^{2}u^{6}}{25b^{2}}
>y_{m_{1}}
>\frac{59.2 \cdot 67.1^{1/(2r_{0}-1)}\left| N_{\alpha} \right|^{3/2+5/(2(2r_{0}-1))} b^{7/(2r_{0}-1)}}
{d^{1/2+1/(2r_{0}-1)} \left( |g|\cN_{d',4} \right)^{2-1/(2r_{0}-1)}},
\]
cannot hold. I.e., we get a contradiction and $r_{0}>1$, $\zeta_{4}p/q \neq p_{r_{0}}/q_{r_{0}}$
for all $4$-th roots of unity, $\zeta_{4}$ is not possible.

Using $|g|\cN_{d',4} \geq 2$ and rearranging terms, this inequality holds if
\begin{equation}
\label{eq:dLB-step3a}
d >\frac{1480^{2(2r_{0}-1)/(10r_{0}-3)} \cdot 67.1^{2/(10r_{0}-3)}
         \left| N_{\alpha} \right|^{2(r_{0}+2)/(10r_{0}-3)}
         b^{2(4r_{0}+5)/(10r_{0}-3)}}
{u^{12(2r_{0}-1)/(10r_{0}-3)} 2^{2(4r_{0}-3)/(10r_{0}-3)}}.
\end{equation}

For $r_{0} \geq 2$, we have the following:\\
$2\left( r_{0}+2 \right)/ \left( 10r_{0}-3 \right)=1/5+23/ \left( 5 \left( 10r_{0}-3 \right) \right)$
decreases from $8/17$ towards $1/5$, so\\
$\left| N_{\alpha} \right|^{2(r_{0}+2)/(10r_{0}-3)} \leq \left| N_{\alpha} \right|^{8/17}$,

\noindent
$2 \left( 4r_{0}+5 \right)/ \left( 10r_{0}-3 \right)=4/5+62/ \left( 5 \left( 10r_{0}-3 \right) \right)$
decreases from $26/17$ towards $4/5$, so\\
$b^{2(4r_{0}+5)/(10r_{0}-3)} \leq b^{26/17}$,

\noindent
$12 \left( 2r_{0}-1 \right)/ \left( 10r_{0}-3 \right)=12/5-24/\left( 5 \left( 10r_{0}-3 \right) \right)$
increases from $36/17$ towards $12/5$,
so $u^{12(2r_{0}-1)/(10r_{0}-3)} \geq u^{36/17}$.

Also, $2 \left( 2r_{0}-1 \right)/ \left( 10r_{0}-3 \right)=2/5-4/\left( 5 \left( 10r_{0}-3 \right) \right)$
and $2 \left( 4r_{0}-3 \right)/ \left( 10r_{0}-3 \right)=4/5-18/\left( 5 \left( 10r_{0}-3 \right) \right)$,
so
\[
\frac{1480^{2(2r_{0}-1)/(10r_{0}-3)} \cdot 67.1^{2/(10r_{0}-3)}}
     {2^{2(4r_{0}-3)/(10r_{0}-3)}}
=\frac{1480^{2/5}}{2^{4/5}}
\left( \frac{67.1^{10} \cdot 2^{18}}{1480^{4}} \right)^{1/(5(10r_{0}-3))}.
\]

Since $67.1^{10} \cdot 2^{18}/1480^{4}>1$, we find that
$\left( 67.1^{10} \cdot 2^{18}/1480^{4} \right)^{1/(5(10r_{0}-3))}$
takes its maximum value at $r_{0}=2$, where it is less than $1.35$.
Since $1480^{2/5}/2^{4/5}<10.65$, we see that \eqref{eq:dLB-step3a} holds if
\begin{equation}
\label{eq:thm13-step-iii-dLB1}
d >\frac{15\left| N_{\alpha} \right|^{8/17} b^{26/17}}
{u^{36/17}}.
\end{equation}

If we apply $\left| N_{\alpha} \right|<db^{4}$ to \eqref{eq:thm13-step-iii-dLB1}, we obtain
\begin{equation}
\label{eq:thm13-step-iii-dLB2}
d>\frac{167 b^{58/9}}{u^{4}}.
\end{equation}

\subsection{$r_{0}>1$ and $\zeta_{4}p/q = p_{r_{0}}/q_{r_{0}}$ for some $4$-th root of unity, $\zeta_{4}$}
\label{subsect:thm13-step-iv}

As in the previous subsection, we want to show that
\[
\frac{|N_{\alpha}|d^{2}u^{6}}{25b^{2}},
\]
the lower bound for $y_{m_{1}}$ from \eqref{eq:ym1-LB}, is larger than the upper
bound for $y_{m_{1}}$ in \eqref{eq:yk-UB-step-iva}. That is, we want to show that
\[
d^{5/2+1/(2(r_{0}-1))}
>65.32\frac{b^{2+2/(r_{0}-1)}\left| N_{\alpha} \right|^{1/2+3/(2(r_{0}-1))} \cdot 422^{1/(2(r_{0}-1))}}
{u^{6}|g|^{2}\cN_{d',4}^{2}}.
\]

This becomes
\[
d>65.32^{2/5-2/(5(5r_{0}-4))}\frac{b^{4/5+16/(5(5r_{0}-4))}
\left| N_{\alpha} \right|^{1/5+14/(5(5r_{0}-4))} \cdot 422^{1/(5r_{0}-4)}}
{u^{12/5-12/(5(5r_{0}-4))}\left( |g|\cN_{d',4} \right)^{4/5-4/(5(5r_{0}-4))}}.
\]

As in Subsection~\ref{subsect:thm13-step-iii}, for $r_{0} \geq 2$,\\
$4/5+16/ \left( 5 \left( 5r_{0}-4 \right) \right) \leq 4/3$, with the maximum value occurring for $r_{0}=2$,\\
$1/5+14/\left( 5 \left( 5r_{0}-4 \right) \right) \leq 2/3$, with the maximum value occurring for $r_{0}=2$,\\
$12/5-12/\left( 5 \left( 5r_{0}-4 \right) \right) \geq 2$, with the minimum value occurring for $r_{0}=2$.

Applying $|g|\cN_{d',4} \geq 2$, we have
\[
\frac{65.32^{2/5}}{2^{4/5}}
\left( \frac{422^{5} \cdot 2^{4}}{65.32^{2}} \right)^{1/(5(5r_{0}-4))}
< 3.057 \cdot 1.472 < 4.5,
\]
for $r_{0} \geq 2$.

Hence
\begin{equation}
\label{eq:thm13-step-iv-dLB1}
d>\frac{4.5b^{4/3}\left| N_{\alpha} \right|^{2/3}}{u^{2}}.
\end{equation}

If we apply $\left| N_{\alpha} \right|<db^{4}$ to \eqref{eq:thm13-step-iv-dLB1}, we obtain
\begin{equation}
\label{eq:thm13-step-iv-dLB2}
d > \frac{92b^{12}}{u^{6}}.
\end{equation}

\subsection{Proof of Theorem~\ref{thm:1.3-seq-new}(b)}
\label{subsect:proof-thm13b}

We now bring together all of the conditions that we have imposed to obtain a
contradiction from the assumption that for $d$ sufficiently large, there are
five distinct squares with $k \geq 3$ or $k \leq K-2$.

From \eqref{eq:ym1-LB}, for such $m_{1}$, we have
$y_{m_{1}} \geq \left| N_{\alpha} \right| d^{2}u^{6} / \left( 25b^{2} \right)$.
For \eqref{eq:yks-LB2}, we need 
$y_{m_{1}} \geq 4\left| N_{\alpha} \right|/\sqrt{d}$. So this holds if
\begin{equation}
\label{eq:dUB-from-ass-2a}
d \geq \left( \frac{10000b^{4}}{u^{12}} \right)^{1/5}.
\end{equation}

We also need $y_{m_{1}} \geq (64/15)b^{4} \left| N_{\alpha} \right|^{4}/d^{2}$ for \eqref{eq:yks-LB2}.
So using \eqref{eq:ym1-LB} again, we need to show that
$\left| N_{\alpha} \right| d^{2}u^{6}/ \left( 25b^{2} \right)
>(64/15)b^{4} \left| N_{\alpha} \right|^{4}/d^{2}$.
This holds if
\begin{equation}
\label{eq:dUB-from-ass-2b-a}
d \geq \frac{(320/3)^{1/4}\left| N_{\alpha} \right|^{3/4}b^{3/2}}{u^{3/2}}.
\end{equation}

As with the lower bounds for $d$ in the previous subsections, it will be useful
to have a lower bound for $d$ here that depends only on $b$ and $u$. Using
$\left| N_{\alpha} \right| < db^{4}$, we find that
\begin{equation}
\label{eq:dUB-from-ass-2b-b}
d \geq \frac{320b^{18}}{3u^{6}}.
\end{equation}

In Subsection~\ref{subsect:thm13-step-i}, we also added the assumption that
$d>8.8b^{390/107}/u^{330/107}$ in \eqref{eq:thm13-step-i-dLB2}.

In Subsection~\ref{subsect:thm13-step-ii}, we also added the assumption that
$d > 6b^{98/25}/u^{78/25}$ in \eqref{eq:thm13-step-ii-dLB2}.

At the end of Subsection~\ref{subsect:thm13-step-iii}, we imposed the condition that
$d>167 b^{52/9}/u^{4}$ in \eqref{eq:thm13-step-iii-dLB2}.

Finally, at the end of Subsection~\ref{subsect:thm13-step-iv}, we imposed the condition
that $d>3775b^{12}/u^{6}$ in \eqref{eq:thm13-step-iv-dLB2}.

Combining these lower bounds for $d$, we have
\begin{equation}
\label{eq:d-LB2}
d \geq \max \left(
\left( \frac{10000b^{4}}{u^{12}} \right)^{1/5},
\frac{320b^{18}}{3u^{6}},
\frac{8.8b^{390/107}}{u^{330/107}}, % (AS-3b)
\frac{6b^{98/25}}{u^{78/25}}, % (AS-4b)
\frac{167 b^{52/9}}{u^{4}}, % (AS-5b)
\frac{92b^{12}}{u^{6}} \right). % (AS-6b)
\end{equation}

If instead of \eqref{eq:thm13-step-i-dLB2}, \eqref{eq:thm13-step-ii-dLB2}, \eqref{eq:thm13-step-iii-dLB2} and
\eqref{eq:thm13-step-iv-dLB2}, we use \eqref{eq:thm13-step-i-dLB1}, \eqref{eq:thm13-step-ii-dLB1}, \eqref{eq:thm13-step-iii-dLB1}
and \eqref{eq:thm13-step-iv-dLB1} and round up some of the constants, then we obtain
\[
d \geq \max \left(
\frac{7b^{4/5}}{u^{12/5}},
\frac{3.22\left| N_{\alpha} \right|^{3/4}b^{3/2}}{u^{3/2}},
\frac{4.6 \left| N_{\alpha} \right|^{45/152} b^{105/76}}{u^{165/76}}, % (AS-3a)
\frac{3.46 \left| N_{\alpha} \right|^{11/36} b^{3/2}}{u^{13/6}}, % (AS-4a)
\frac{15 \left| N_{\alpha} \right|^{8/17}b^{20/17}}{u^{36/17}}, % (AS-6a)
\frac{4.5 \left| N_{\alpha} \right|^{2/3}b^{4/3}}{u^{2}}
\right).
\]

Notice that the exponents on both $b$ and $\left| N_{\alpha} \right|$ in the
second term in the max are at least as large as those on the other terms.
Also the exponent on $u$ in the second term is smaller than in the other terms
and the largest constant on the terms is $15$. Hence
\[
d \geq
\frac{15\left| N_{\alpha} \right|^{3/4} b^{3/2}}{u^{3/2}}.
\]

This is the lower bound for $d$ in Theorem~\ref{thm:1.3-seq-new}(b), so this
completes the proof of that part of the theorem.

\subsection{Proof of Theorem~\ref{thm:1.3-seq-new}(a)}
\label{subsect:proof-thm13a}

As in Subsection~\ref{subsect:proof-thm13b}, we are going to obtain a contradiction
from the assumption that there are five distinct squares, $y_{k}$, with $k \geq 3$
or $k \leq K-2$.

To complete the proof of Theorem~\ref{thm:1.3-seq-new}(a), we need to remove the
assumption on $d$ in \eqref{eq:d-LB2} in Subsection~\ref{subsect:proof-thm13b} for
$b=1,2,3$. To make the calculation for $b=3$ more feasible, we will refine
\eqref{eq:d-LB2} somewhat by improving \eqref{eq:dUB-from-ass-2b-b}.

\begin{lemma}
\label{lem:improved-dLB}
If $d>201$ and $d>9.47b^{18}/u^{6}$, then $y_{m_{1}} \geq (64/15)b^{4} \left| N_{\alpha} \right|^{4}/d^{2}$.
\end{lemma}

\begin{proof}
Rather than using \eqref{eq:ym1-LB} as in Subsection~\ref{subsect:proof-thm13b},
we will use Lemma~\ref{lem:Y-LB1}(a) instead. Since $0<-\overline{\alpha}
=b^{2}\sqrt{d}-a < \alpha=a+b^{2}\sqrt{d}$, we want to show that
\begin{equation}
\label{eq:init-ineq}
\frac{-\overline{\alpha} \varepsilon^{6}}{2\sqrt{d}}
>\frac{(64/15)b^{4} \left| N_{\alpha} \right|^{4}}{d^{2}}
\end{equation}
holds when $d>9.47b^{18}/u^{6}$.

This inequality holds when
\begin{equation}
\label{eq:lemma-dLB1}
d>\frac{4.176b^{8/3} \left| N_{\alpha} \right|^{8/3}}{(-\overline{\alpha})^{2/3} \varepsilon^{4}}.
\end{equation}

For $d>201$, we have $t^{2}\geq du^{2}-4 \geq 0.9801du^{2}$. So
$\varepsilon>0.995\sqrt{d} \, u$. Substituting this lower bound for $\varepsilon$
into \eqref{eq:lemma-dLB1} and then solving for $d$, we see that the desired
inequality holds for
\begin{equation}
\label{eq:lemma-dLB2}
d>\frac{1.622b^{8/9} \left| N_{\alpha} \right|^{8/9}}{(-\overline{\alpha})^{2/9} u^{4/3}}
=\frac{1.622b^{8/9} \left| N_{\alpha} \right|^{2/3} \alpha^{2/9}}{u^{4/3}}.
\end{equation}

Next, taking the derivative of $\left| N_{\alpha} \right|^{2/3} \alpha^{2/9}
= \left( b^{4}d-a^{2} \right)^{2/3} \left( a+b\sqrt{d} \right)^{2/9}$ with respect to
$a$ and solving for where it is zero, we want to solve
$7a^{2}+6ab^{2}\sqrt{d} - b^{4}d=\left( a+b^{2}\sqrt{d} \right) \left( 7a-b^{2}\sqrt{d} \right)=0$.
For positive $a$, this occurs at $a=b^{2}\sqrt{d}/7$. Substituting this into
the last expression in \eqref{eq:lemma-dLB2}, we find that \eqref{eq:init-ineq}
holds for $d>9.47b^{18}/u^{6}$, proving the lemma.
\end{proof}

Applying this lemma, we will replace \eqref{eq:d-LB2} with
\begin{equation}
\label{eq:d-LB3}
d \geq \max \left(
\left( \frac{10000b^{4}}{u^{12}} \right)^{1/5},
\frac{9.47b^{18}}{u^{6}},
\frac{8.8b^{390/107}}{u^{330/107}}, % (AS-3b)
\frac{6b^{98/25}}{u^{78/25}}, % (AS-4b)
\frac{167 b^{52/9}}{u^{4}}, % (AS-5b)
\frac{92b^{12}}{u^{6}} \right), % (AS-6b)
\end{equation}
for $d>201$.

This reduces the lower bound for $d$ when $b=3$ and $u=1$ by a factor of over $11$.

We use \eqref{eq:d-LB3} to get initial upper bounds, $D_{b}$, on the values of
$d$ that we need to check.

\begin{center}
\begin{table}[h]
\begin{tabular}{|c|r|r|r|r|}\hline
$b$ &     $D_{b}$          &       $c_{1,b}$ & $c_{2,b}$ & \text{CPU time} \\ \hline
 1  &                $167$ &           $144$ &     $2$   & $<1$s \\ \hline
 2  &          $2,482,503$ &     $1,881,618$ &     $0$   &   33s \\ \hline
 3  & $3.669 \cdot 10^{9}$ & $6,238,802,661$ &     $0$   & 64h21m \\ \hline
\end{tabular}
\caption{Data for $b$}
\label{table:b-data}
\end{table}
\end{center}

\begin{center}
\begin{table}[h]
\begin{tabular}{|c|c|c|}\hline
  $(a,b,d,t,u)$    & $N_{\alpha}$ \\ \hline
$(2, 1,  8, 2, 1)$ &      $-4$    \\ \hline
$(2, 1, 13, 3, 1)$ &      $-9$    \\ \hline
\end{tabular}
\caption{Remaining tuples}
\label{table:eqn-data}
\end{table}
\end{center}

For $b=1,2$, the value of $D_{b}$ in Table~\ref{table:b-data} comes
from the sixth term in the max in \eqref{eq:d-LB3} with $u=1$.
For $b=3$, the values of $b$, the value of $D_{b}$ in Table~\ref{table:b-data}
comes from the second term in the max in \eqref{eq:d-LB3} with $u=1$. It
is over twice the size of the value of the sixth term in the max in \eqref{eq:d-LB3}
with $u=1$. Note that for all three, the lower bound for $d$ exceeds the hypothesis
in Lemma~\ref{lem:improved-dLB}, so the use of the \eqref{eq:d-LB3} is permitted
here.

For each value of $b$, find the smallest value of $u$ such that the right-hand
side of \eqref{eq:d-LB2} is less than $2$ and denote this value of $u$ by $U_{b}$.

For each $u$ with $1 \leq u \leq U_{b}$, consider all integers
$d=\left( t^{2} \pm 4 \right)/ u^{2} \leq D_{b}$.

For each such $(b,d,t,u)$, we find all positive integers, $a$, such that $N_{\alpha}<0$
and such that \eqref{eq:dUB-from-ass-2b-a}, \eqref{eq:thm13-step-i-dLB1},
\eqref{eq:thm13-step-ii-dLB1}, \eqref{eq:thm13-step-iii-dLB1} and \eqref{eq:thm13-step-iv-dLB1}
all hold. In fact, only \eqref{eq:dUB-from-ass-2b-a} yields a non-trivial bound.
The others give a non-positive lower bound for $a$. From \eqref{eq:dUB-from-ass-2b-a},
we obtain
$db^{4}-a^{2}=\left| N_{\alpha} \right| \leq \left( 3d^{4}u^{6}/ \left( 320b^{6} \right) \right)^{1/3}$
and hence $a \geq \sqrt{db^{4}-\left( 3d^{4}u^{6}/ \left( 320b^{6} \right) \right)^{1/3}}$.
We also have the upper bound $a<\sqrt{db^{4}}$ from
$N_{\alpha}<0$. We found $c_{1,b}$ such tuples $(a,b,d,t,u)$, where $c_{1,b}$ is recorded
in Table~\ref{table:b-data}.

For each of these tuples $(a,b,d,t,u)$, we checked for squares among the $y_{k}$'s
satisfying \eqref{eq:yks-LB2}, along with the upper bounds in \eqref{eq:yk-UB-step-i},
\eqref{eq:yk-UB-step-ii}, \eqref{eq:yk-UB-step-iii} and \eqref{eq:yk-UB-step-iv}.

If such a square, $y_{k}$, exists, then the tuple $(a,b,d,t,u)$ requires further
work that we conduct separately in Magma, using its \verb+IntegralQuarticPoints()+
function. Of the above $c_{1,b}$ tuples $(a,b,d,t,u)$,
all but $c_{2,b}$ tuples (where $c_{2,b}$ is as in Table~\ref{table:b-data}) are
filtered out using these checks. The remaining two tuples are given in
Table~\ref{table:eqn-data}. Since $-N_{\alpha}$ is a positive square for both
of these tuples, by Theorem~\ref{thm:any-b-sqr}(b), part~(a) of the theorem
holds for these tuples.
Hence part~(a) of the theorem holds.

Code for this search for such tuples $(a,b,d,t,u)$, was written
in PARI/GP and run on a Windows laptop with an Intel i7-13700H CPU and 32GB of
memory. We record the CPU time for each value of $b$ in Table~\ref{table:b-data}.
The PARI/GP code can be found in the file \verb!d-bAll-nAAny-search.gp!
in the pari subdirectory of the github url provided at the end of Section~\ref{sect:intro}.


\begin{thebibliography}{10}
\bibitem{Pari}
The PARI Group, PARI/GP version {\tt 2.16.2}, Univ. Bordeaux, 2024,
\url{http://pari.math.u-bordeaux.fr/}.

\bibitem{V5}
P. M. Voutier,
\emph{Bounds on the number of squares in recurrence sequences},
J. Number Theory {\bf 265} (2024), 291--343.

\bibitem{V6}
P. M. Voutier,
\emph{Sharp bounds on the number of squares in recurrence sequences and solutions of $X^{2}-\left( a^{2}+b \right) Y^{4}=-b$},
Research in Number Theory (accepted) \url{https://arxiv.org/abs/1807.04116}.

\bibitem{V7}
P. M. Voutier,
\emph{Bounds on the number of squares in recurrence sequences: arbitrary $b$, I},
(submitted) \url{http://arxiv.org/abs/2502.14875}.

\bibitem{V8}
P. M. Voutier,
\emph{Bounds on the number of squares in recurrence sequences: arbitrary $b$, II},
(submitted).
\end{thebibliography}
\end{document}